\documentclass[11pt]{article}
\usepackage[utf8]{inputenc}
\usepackage[T1]{fontenc}
\usepackage{lmodern}
\usepackage{amsmath}
\usepackage{amssymb}
\usepackage{amsthm}
\usepackage{booktabs}
\usepackage{array}
\usepackage{geometry}
\usepackage{microtype}
\usepackage{hyperref}
\usepackage{url}
\usepackage{enumitem}
\usepackage{seqsplit}
\hypersetup{colorlinks=true,linkcolor=blue,	citecolor=blue,	urlcolor=blue}
\title{\textbf{When Fibonacci and Lucas Meet Smith:\\
A Computational Exploration}}
\author{Ronie Peterson Dario	\qquad	Jo\~ao Luis Gon\c{c}alves \qquad Moniky P. N. De Oliveira
\\[0.6em] 	\small Universidade Tecnol\'ogica Federal do Paran\'a (UTFPR) \\ 	\small Curitiba, PR, Brazil
}
\date{}
\begin{document}
\maketitle
\begin{abstract}The Fibonacci and Lucas sequences are old friends in recreational number theory. Here we ask a simple question: which of their terms \(F_n\) and \(L_n\) are Smith numbers? A Smith number is a composite integer whose decimal digit sum equals the sum of the decimal digits of its prime factors, counted with multiplicity. The question is easy to state, but for large \(n\) it quickly becomes a factorization problem.

Using available complete factorization data, we carried out computational searches in both sequences and obtained new Smith terms in each of them.

For the Lucas sequence, our search found ten Smith terms with indices below \(1000\):
\[
3,\ 95,\ 105,\ 114,\ 183,\ 437,\ 609,\ 682,\ 827,\ 902.
\]
These results were published in the On-Line Encyclopedia of Integer Sequences as OEIS A395686. Four further indices,
\[1090,\ 1153,\ 1215,\ 1378,\]
were subsequently added to the sequence by Sean A. Irvine.

For the Fibonacci sequence, starting from the cases already recorded in OEIS A382922, our search produced seven further Smith numbers:
\[
F_{1440},\quad
F_{1554},\quad
F_{1596},\quad
F_{1863},\quad
F_{2256},\quad
F_{2277},\quad
F_{2559}.
\]

For the even-index cases, the identity $F_{2n}=F_nL_n$ allows available Fibonacci and Lucas factorization data to be combined in order to recover the complete factorizations needed for the Smith test. Surprisingly, the two searches also meet at the index \(827\): both
\[F_{827} \text{ \ and \ } L_{827}\] are Smith numbers.
\end{abstract}

\section{A telephone number meets Fibonacci and Lucas}

Smith numbers entered mathematics through an appropriately
recreational route. In 1982 Albert Wilansky introduced a curious class of integers inspired by the telephone number of his brother-in-law,	Harold Smith \cite{wilansky}. The number was
\[	4937775.	\]
Its prime factorization is \[4937775=3\cdot5\cdot5\cdot65837.\]
Adding its decimal digits gives
\[4+9+3+7+7+7+5=42,	\]
while doing the same with the prime factors, counted with multiplicity, gives
\[	3+5+5+(6+5+8+3+7)=42.	\]

This coincidence gave rise to the name \emph{Smith number}.
The definition combines two different arithmetic features of an
integer: its decimal representation and its prime factorization. It is therefore natural to ask how the Smith property behaves within familiar integer sequences.

Consider the Fibonacci sequence
\begin{equation}\label{eq:fibo-sequence-defn}	F_0=0,\qquad
F_1=1,\qquad 	F_{n}=F_{n-1}+F_{n-2},\ \ n\geq 2,
\end{equation}
and its close relative, the Lucas sequence
\begin{equation}\label{eq:lucas-sequence-defn}
L_0=2,\qquad 	L_1=1,\qquad 	L_{n}=L_{n-1}+L_{n-2}, \ \ n\geq 2.
\end{equation}

Large Fibonacci and Lucas numbers are easily generated, but testing the Smith property hides a much harder computational problem. Adding the decimal digits of a term with hundreds of digits is straightforward; obtaining its complete prime factorization may be a substantial computation, and
in many cases the factorization is not yet known.

This leads to the question considered throughout this paper:
\begin{center}
\emph{Which Fibonacci and Lucas numbers are Smith numbers, and how far
can this question be explored computationally?}
\end{center}

We first recall the Smith condition and then examine the Lucas and Fibonacci sequences separately. The classical identity
\begin{equation}\label{eq:F_2n=F_nL_n} F_{2n}=F_nL_n \end{equation}
will later provide a useful connection between the two searches.

\section{The Smith condition}\label{sec:smith-condition}

For an integer \(N\geq 1\), denote by \(S(N)\) the sum of its decimal
digits. If
\[
N=p_1^{a_1}p_2^{a_2}\cdots p_r^{a_r}>1
\]
is its prime factorization, set
\[
P(N)=\sum_{i=1}^{r}a_iS(p_i).
\]
Thus \(P(N)\) is the sum of the decimal digit sums of the prime factors
of \(N\), counted with multiplicity.

A composite integer \(N\) is a \emph{Smith number} if
\begin{equation}\label{eq:smith-condition-defn} S(N)=P(N).\end{equation}

For example, $666=2\cdot3^2\cdot37$ is a Smith number, since
\[S(666)=18 \text{ \ and \ } S(2)+2S(3)+S(37) = 2+6+10 = 18.\]

Prime numbers are excluded from the definition, since for a prime \(p\) the equality \eqref{eq:smith-condition-defn} would hold automatically.

A useful feature of $P$ is its additivity under multiplication:
\begin{equation}\label{eq:additive-of-P}
P(ab)=P(a)+P(b) \end{equation}
for integers \(a\), \(b > 1\). This simple observation will later be combined with \eqref{eq:F_2n=F_nL_n}.

Smith numbers form an infinite set, as proved by McDaniel
\cite{mcdaniel}. Our concern here is their occurrence within the
Fibonacci and Lucas sequences.

\section{The Lucas--Smith search}\label{sec:lucas-smith-search}

We first consider the Lucas sequence. Using the Fibonacci and Lucas factorization tables maintained by
B. Kelly \cite{kelly}, together with the classical tables of
Brillhart, Montgomery, and Silverman \cite{bms}, we examined those terms \(L_n\) for which a complete prime factorization was available or could be reconstructed from the listed data.

For \(n<1000\), the search identified ten indices for which \(L_n\) is a Smith number:
\[3,\ 95,\ 105,\ 114,\ 183,\ 437,\ 609,\ 682,\ 827,\ 902.\]

These results were published in the On-Line Encyclopedia of Integer Sequences as OEIS A395686 \cite{oeis-lucas-smith}. Sean A. Irvine subsequently extended the sequence by adding four further indices: $1090$, $1153$, $1215$, $1378.$
In the factorization data used here, complete factorizations are available through \(n=1422\). The case \(n=1423\) is the first for which the available factorization is incomplete, and therefore provides a natural stopping point for this part of the search.

To illustrate the computation behind these data, consider \(L_{437}\). It is the \(92\)-digit integer
\[
\begin{aligned}
L_{437}={}&
2126173932344343953551107295983794236927231713\\
&3885364875092188438605064381241873859092393771.
\end{aligned}
\]
with prime factorization
\[\begin{aligned}
L_{437}={}&139\cdot461\cdot9349 \cdot133351074933356474419349\\
&\cdot266146882200264269210394982818118404450063235825822769857949.
\end{aligned}\]
Its decimal digit sum is $S(L_{437})=416,$ while
\[P(L_{437}) = S(139)+S(461)+S(9349) +S(1...9) + S(2...9) = 13+11+25+101+266 = 416.\]
Thus \(S(L_{437})=P(L_{437})\), confirming that \(L_{437}\) is a Smith number.

Table~\ref{tab:lucas-smith} summarizes the numerical data for the
fourteen Smith Lucas numbers currently recorded in this range.

\begin{table}[htbp]
\centering
\caption{Smith numbers found in the Lucas sequence.}
\label{tab:lucas-smith}

\begin{tabular}{rcc}
\toprule
\(n\) &
Digits of \(L_n\) &
\(S(L_n)=P(L_n)\)
\\
\midrule
3    &   1 &    4\\
95   &  20 &   89\\
105  &  22 &   76\\
114  &  24 &  108\\
183  &  39 &  194\\
437  &  92 &  416\\
609  & 128 &  652\\
682  & 143 &  618\\
827  & 173 &  802\\
902  & 189 &  816\\
1090 & 228 &  987\\
1153 & 241 & 1117\\
1215 & 254 & 1148\\
1378 & 288 & 1293\\
\bottomrule
\end{tabular}
\end{table}

\section{Smith numbers hiding among Fibonacci numbers}
\label{sec:smith-numbers-in-fibo-sequence}

We now turn to the Fibonacci sequence. A relatively small example is
\[F_{77}=5527939700884757,\]
with prime factorization
\[F_{77} = 13\cdot89\cdot988681\cdot4832521.\]
Its decimal digit sum is $S(F_{77})=86,$ and the digit sums of its prime factors also add to \(86\). Hence \(F_{77}\) is a Smith number.

The indices \(n\) for which \(F_n\) is a Smith number are recorded in OEIS A382922 \cite{oeis-fibonacci-smith}. Before the computational extension considered here, the listed indices were
\[ 31,\ 77,\ 231,\ 354,\ 523,\ 535,\ 631,\ 819,\ 827,\ 830,\ 991,\ 1234. \]

For the Fibonacci search, we again relied on the factorization data  available in B. Kelly's tables \cite{kelly} and in the tables of Brillhart, Montgomery, and Silverman \cite{bms}. Whenever the complete prime factorization of \(F_n\) was available or could be reconstructed from these data, we compared \(S(F_n)\) and \(P(F_n)\).

For odd indices, the necessary factors can often be recovered directly from the Fibonacci factorization data. For even indices, Lucas factorization data also enter through the identity \eqref{eq:F_2n=F_nL_n}; this connection will be discussed in the next section.

In the range
\[1235\leq n<10\,000,\]
the search produced seven further Smith Fibonacci numbers:
\[
F_{1440},\quad F_{1554},\quad F_{1596},\quad F_{1863},\quad F_{2256},\quad  F_{2277},\quad F_{2559}.\]

Their numerical data are summarized in Table~\ref{tab:fibonacci-smith}.

\begin{table}[htbp]
\centering
\caption{Further Smith numbers found in the Fibonacci sequence.}
\label{tab:fibonacci-smith}

\begin{tabular}{rcc}
\toprule
\(n\) &
Digits of \(F_n\) &
\(S(F_n)=P(F_n)\)
\\
\midrule
1440 & 301 & 1296\\
1554 & 325 & 1378\\
1596 & 334 & 1557\\
1863 & 389 & 1780\\
2256 & 472 & 2115\\
2277 & 476 & 2144\\
2559 & 535 & 2473\\
\bottomrule
\end{tabular}
\end{table}

To illustrate the scale of the computation, consider \(F_{1863}\). It has \(389\) decimal digits, and its complete prime factorization can be written as
\begingroup
\small
\begin{align*}
F_{1863}={}& 2\cdot17\cdot53\cdot109\cdot137\cdot829
\cdot2269\cdot4373\cdot18077\cdot19441
\cdot28657\cdot33533\cdot10656361
\cdot23019229 \\ & \cdot840175741  \cdot17642247580301401 \cdot1497034599690690354913433 \\ &  \cdot4072353155773627601222196481\\
&\cdot 305504674058681335357289317796420951893158199699949280927361 \\ & \cdot \mathrm{P}204\end{align*}
\endgroup
where $\mathrm{P}204$ denotes a \(204\)-digit prime.  This factor can be recovered by generating \(F_{1863}\) from the recurrence \eqref{eq:fibo-sequence-defn} and dividing by the product of the
other prime factors displayed above. Despite the size of the factorization, the final Smith test is simply
\[S(F_{1863})=P(F_{1863})=1780. \]

\section{The link between the two searches}

The two searches are connected by the classical identity \eqref{eq:F_2n=F_nL_n}.  For even Fibonacci indices, this allows Fibonacci and Lucas
factorization data to be combined. Since \(P\) is additive under
multiplication,
\[P(F_{2m})=P(F_m)+P(L_m).\]

The case \(n=1440\) provides a simple illustration. Since
\[F_{1440}=F_{720}L_{720},\]
the available complete factorizations give
\[
P(F_{720})=615
\qquad\text{and}\qquad
P(L_{720})=681.
\]
Hence
\[
P(F_{1440})=615+681=1296.
\]
On the other hand, $S(F_{1440})=1296$. Hence \(F_{1440}\) is a Smith number.

The same procedure applies to the other even-index Fibonacci terms
identified in the search.

\section{A surprising common index}\label{sec:surprise}

Among the cases identified in the two searches, one index appears in both lists: $827$. Indeed,
\[F_{827} \text{ \ \ and  \ \ } L_{827}\] are both Smith numbers.

Thus the Fibonacci and Lucas searches meet at the same index in a
rather unexpected way. This naturally raises a further question:

\begin{center}
\emph{Are there other indices \(n\) for which both \(F_n\) and
\(L_n\) are Smith numbers?}
\end{center}

\section{Final remarks}

We began with a simple question: which Fibonacci and Lucas numbers are Smith numbers, and how far can this question be explored computationally? The answer quickly led from decimal digit sums to the much harder problem of complete prime factorization.

For the Lucas sequence, the search produced ten Smith terms below index \(1000\), leading to OEIS A395686. Sean A. Irvine subsequently added four further terms. For the Fibonacci sequence, our search produced seven further Smith numbers beyond those previously recorded in OEIS A382922.

The identity \eqref{eq:F_2n=F_nL_n} provides a useful bridge between the two searches, allowing Lucas factorization data to contribute directly to the analysis of even-index Fibonacci numbers. The most unexpected meeting, however, occurs at \(n=827\): both \(F_{827}\) and \(L_{827}\) are Smith numbers. Whether Fibonacci and Lucas meet Smith again at the same
index remains an open question.

There is also room for both lists to grow. Many large Fibonacci and Lucas numbers are not yet completely factored, and therefore cannot yet be definitively tested for the Smith property. A new complete factorization may consequently reveal another Smith number hidden in one of these familiar sequences.

We are especially grateful to B. Kelly for maintaining and extending the Fibonacci and Lucas factorization tables, and sincerely hope that this valuable computational effort will continue. We also thank Sean A. Irvine for valuable discussions and for his interest in the Lucas--Smith sequence.


\begin{thebibliography}{BMS88}

\bibitem[BMS88]{bms}
J. Brillhart, P. L. Montgomery, and R. D. Silverman.
Tables of Fibonacci and Lucas factorizations.
\emph{Mathematics of Computation}, 50:251--260, 1988.

\bibitem[Cos02]{costello}
P. Costello.
A new largest Smith number.
\emph{The Fibonacci Quarterly}, 40(4):369--371, 2002.

\bibitem[Cos15]{costello2015}
P. Costello.
Smith numbers from primes with small digits.
\emph{Missouri Journal of Mathematical Sciences},
27(1):10--15, 2015.

\bibitem[Kel]{kelly}
B. Kelly.
Fibonacci and Lucas factorizations.
\url{https://mersennus.net/fibonacci/}.
Accessed August 23, 2026.

\bibitem[Kos01]{koshy}
T. Koshy.
\emph{Fibonacci and Lucas Numbers with Applications}.
Wiley, New York, 2001.

\bibitem[McD87]{mcdaniel}
W. L. McDaniel.
The existence of infinitely many \(k\)-Smith numbers.
\emph{The Fibonacci Quarterly}, 25(1):76--80, 1987.

\bibitem[OEIS1]{oeis-smith}
N. J. A. Sloane and The OEIS Foundation Inc.
Smith (or joke) numbers.
\emph{The On-Line Encyclopedia of Integer Sequences},
A006753.
\url{https://oeis.org/A006753}.

\bibitem[OEIS2]{oeis-fibonacci-smith}
S. S. Gupta.
Numbers \(k\) such that Fibonacci(\(k\)) is a Smith number.
\emph{The On-Line Encyclopedia of Integer Sequences},
A382922.
\url{https://oeis.org/A382922}.

\bibitem[OEIS3]{oeis-lucas-smith}
R. P. Dario.
Indices \(n\) for which the Lucas number \(L_n\) is a Smith number.
\emph{The On-Line Encyclopedia of Integer Sequences},
A395686.
\url{https://oeis.org/A395686}.

\bibitem[Wil82]{wilansky}
A. Wilansky.
Smith numbers.
\emph{Two-Year College Mathematics Journal}, 13:21, 1982.

\end{thebibliography}
\end{document}